\documentclass[twocolumn,pre]{revtex4}

\usepackage{dcolumn}
\usepackage{amsmath}

\usepackage{graphicx}

\newlength{\figurewidth}
\begin{document}

\title{The structure of self-organized blogosphere}
\author{Feng Fu$^1$}
\email{fufeng@pku.edu.cn}
\author{Lianghuan Liu$^1$}
\author{Kai Yang$^2$}
\author{Long Wang$^1$}
\email{longwang@pku.edu.cn} \affiliation{ $^1$Laboratory for
Intelligent Control, Center for Systems and Control, College of
Engineering,\\ Peking University, Beijing 100871, P.R. China\\
$^2$Department of Computer Science and Technology, Tsinghua
University, Beijing 100084, P.R. China }
\begin{abstract}
In this paper, a statistical analysis of the structure of one blog
community, a kind of social networks, is presented. The quantities
such as degree distribution, clustering coefficient, average
shortest path length are calculated to capture the features of the
blogging network. We demonstrate that the blogging network has
small-world property and the in and out degree distributions have
power-law forms. The analysis also confirms that blogging networks
show in general disassortative mixing pattern. Furthermore, the
popularity of the blogs is investigated to have a Zipf's law,
namely, the fraction of the number of page views of blogs follows
a power law.
\end{abstract}

\maketitle

\subsection{Introduction}
The recent development of so-called network science reveals the
underlying structures of complex networks and becomes the catalyst
for arising common voice of interdisciplinary fields to tame the
complexity
\cite{strotagz_nature,Newman,Albert_Barabasi_review_2002,physrep_06}.
The small-world network model proposed by Watts and Strogatz
quantitatively reflected that the real networks are small worlds
which have high clustering and short average path length
\cite{Watts-Strogatz}. The six degree of separation, uncovered by
the social psychologist Stanley Milgram, is the most famous
manifestation of small-world theory \cite{milgram}. The real
world, however, significantly deviated from classic
Erd\"{o}s-R\'{e}nyi model that the degree distribution is right
skew, namely, follows a power law other than Poisson distribution
\cite{ER_1,ER_2,BA_science_1999}. In particular, for most
networks, including the World Wide Web, the Internet and the
metabolic networks, the degree distribution has a power-law tail
--- $p(k)\sim k^{-\gamma}$. Such networks are called scale free and
the Barab\'{a}si-Albert model (BA model) provides a possible
generating mechanism for such scale-free structure: growth and
preferential attachment \cite{BA_science_1999}. These pioneering
discoveries intensively attracted a large number of scientist from
different background to plunge into this emerging immature realm.
Besides, the real networks are hierarchical and have communities
structure or composed of the elements --- motifs
\cite{newman_2006_pnas,motif}. Nevertheless, surprisingly, it is
found that the complex networks are self-similar, corresponding to
the ubiquitous geometry pattern in snowflakes
\cite{self-similarity_nature}. Meanwhile, the dynamics taking
placing on complex networks such as virus spreading (information
propagation), synchronization processes, games and cooperation,
have been deeply investigated and well understood
\cite{Watts-book-sm,info_prop,sf_threshod,kuramoto_rmp,sync_sm,chou_1,chou_2,sf_cooperation}.

The word blog is short for neologism ``Web log'', which is often a
personal journal maintained on the Web. In the past few years,
blogs are the fastest growing part of the WWW \cite{blogistan}.
There are now about 20 million blogs, which are emerging as an
important communication mechanism by an increasing number of
people \cite{nature_blog}. The Web in its first decade was like a
big online library. Today, however, it becomes more of social web,
not unlike Berners-Lee's original vision. Consequently, advanced
social technologies---Blog, Wiki, Podcasting, RSS, etc, which
featured as characteristics of time of Web 2.0 --- have led to the
change of the ways of people's thinking and communicating. We
refer blogistan as blog space in the jargon of the blog field. As
one surfs in blogistan, the global blogistan is just like an
ecosystem called blogosphere that has a life of its own. In the
view of complex adaptive system, the whole blogosphere is more
than the sum of its weblogs. Therefore, one can't understand the
blogosphere by studying one single weblog. Moreover, some
interesting phenomena corresponding to the classic ecological
patterns---predators and prey, evolution and emergence, natural
selection and adaptation---are ubiquitous in blogosphere, where
evolutionary forces plays out in real time. For instance,
individual weblogs vie for niche status, establish communities of
like-minded sites, and jostle links to their sites
\footnote{Available at http://www.microcontentnews.com/articles/\\
blogosphere.htm}. Besides, the fascinating and powerful filtering
effect, namely, collaborative filtering is created by the dynamic
hierarchy of links and recommendations generated by blogs. The
more bloggers there are in a particular community, the more
efficient this filtering becomes, so, counter-intuitively,
reducing information overload \cite{nature_blog}.

A typical blog is one long Web page on content hosting site that
provides blog space. It is basically a large queue with additions
appearing at the top of the page and older material scrolling
down, often partitioned into archives and with links to other
blogs within the same host site (internal links) or to URLs in the
Web (external links). Sometimes, personal blogs could cite
paragraphs of other blogs, often embedded with links that could be
collected by the blog hosting sites and return feedback to the
original bloggers (the term trackback is used in the blog
commnunity). At first glance, blogs are apparently nothing more
than common Web pages. Nevertheless, active blogs are updated with
a frequency significantly higher than a traditional Web page,
often in a bursty manner. The number and quality of links from a
blog are quite different from ordinary Web pages. The links are
updated more frequently by the bloggers and a significant fraction
of the links are to other blogs. Furthermore, blogosphere creates
an instant online communities of diverse topics for bloggers and
readers who could publish their comments on blogs. Therefore it is
more interactive and open than common Web pages. In this sense,
the blogosphere is worth scrutinizing to reveal underlying
mechanism for these interesting phenomena.

In this paper, we concentrate on the sub-ecosystem of global
blogosphere: the blogs hosted by Sina which is the largest Chinese
blog space provider and has about 2 million registered users in
mainland of China \footnote{http://blog.sina.com.cn}. We are
interested in the emerging links pattern between Sina blogs, i.e.,
the collections of links to bloggers' favorite blog sites. For
simplicity, the links out of the domain (http://blog.sina.com.cn)
are omitted. And also, the Zipf's law in popularity of the blogs
is investigated.

The remainder of this paper is organized as follows. Sec.~II deals
with the method of data collection, and Sec.~III performs the
statistical analysis of the structures of such self-organized
blogosphere, including average degree, degree distribution,
clustering coefficient, etc. Finally, Sec.~IV lays out the
conclusion and future work is presented.

\subsection{Data gathering}
Since there are around $20$ million blogs, we focused our eyesight
in a sub-community of global blogosphere---the Chinese blogs
hosted on Sina. We wanted to examine the structures of such
self-organized ``ecosystem'', including the emerging
interconnected pattern of Sina blogs and the Zipf's law in
popularity of blogs. This would be the first stride to explore the
mysteries of vivid blogosphere.

The blog sites of Sina are very regular, and the entry to blog has
two equivalent forms: (a) http://blog.sina.com.cn/m/XXXX, where
XXXX is a string consisted of letters and numbers; (b)
http://blog.sina.com.cn/u/xxxxx, where xxxxx is a $10$ digits
number as user's id. For all users, they have (b) site forms of
their blogs. While for advanced users, they both have (a) and (b)
forms of entries to their blogs. As some bloggers' sites both have
(a) and (b), the mapping relationships between (a) and (b) are
established to avoid the reduplicate results. We designed a simple
WWW robot which began with the most popular blog, which ranked
first in global blogosphere by Technorati
\footnote{http://www.technorati.com}. This popular blog's number
of page views has been more than $30$ million. Along with the
collections of links to favorite Sina blogs, the robot crawled
down a connected networks of $200,399$ nodes, using breadth-first
search method. At the same time, the page views of each visited
blogs were recorded down. Based upon these data, the analysis of
the structure of self-organized blogosphere was carried out in
next section.

\subsection{Results}
\begin{figure}
\includegraphics[width=8cm]{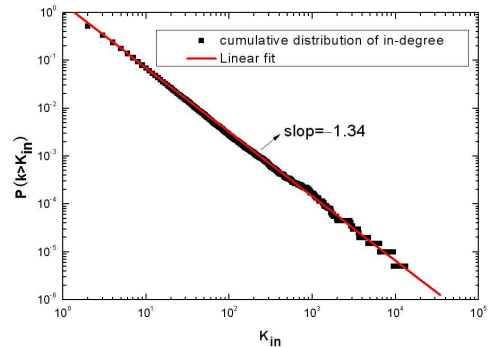}
\caption{\label{fig:indeg}(Color online) Cumulative distribution
of in-degree of blogs. The straight red line is the linear fit of
the data, whose slope is $-1.34\pm 0.001$.}
\end{figure}

\begin{figure}
\includegraphics[width=8cm]{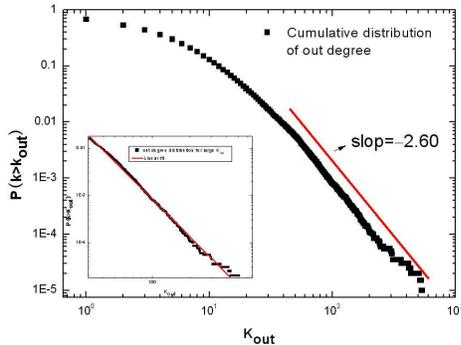}
\caption{\label{fig:outdeg}(Color online) Cumulative distribution
of out degree of blogs. The plotted straight red line is of slope
$-2.60\pm 0.02$ for comparison. The insect shows the detail
log-log plot of the right skew tail for large out degree. The
linear fit of that data justifies that the heavy tail obeys a
power law as $P_{out}(k>K)\sim k^{-\beta}$. }
\end{figure}

The emerging link pattern of Sina blogs is mined from the crawled
down networks. Since the nature of the network is directed, thus
the connectivity of the blog has incoming and outgoing
connections, namely, $k_{in}$ and $k_{out}$ respectively. The
in-degree could be used as an index of importance of the blogs.
From Fig.~\ref{fig:indeg}, we found that the cumulative
distribution of in-degree obeys a power-law form, $P_{in}(k>K)\sim
K^{-\alpha}$, where $\alpha=1.34\pm 0.001$. Therefore, the
in-degree distribution, which indicates the probability that
randomly chosen node $i$ has $k$ incoming connections, follows a
power law, $P_{in}(k)\sim k^{-\gamma_{in}}$, where
$\gamma_{in}=\alpha+1=2.34\pm 0.001$ \footnote{Normally, the right
heavy tail of the distribution is very noisy. To avoid such
difficulty, the cumulative distribution is adopted to measure the
power law correctly. See Ref.~\onlinecite{pl_newman} and
references therein.}. The cumulative distribution of out-degree
has a power-law tail as $P_{out}(k>K)\sim K^{-\beta}$, where
$\beta=2.60\pm 0.02$ (see Fig.~\ref{fig:outdeg} for details).
Thereby, the out-degree distribution has the form $P_{out}(k)\sim
k^{\gamma_{out}}$, where $\gamma=\beta+1=3.60 \pm 0.02$. By
contrast, the out-degree distribution slightly deviates from the
right skew heavy tail for small out-degree $k_{out}$.
Paradoxically, someone may argue that the log-normal distribution
would better fit the data than the power law. Yet, we think that
there exists a threshold as certain $k_{out}$ and when out-degree
exceeds that threshold, a power-law tail exists, given the
evidence that most of the data fall into the right skew tail
\footnote{An alternative explanation is that there exist two
different power-law regimes: one is for sufficiently small
$k_{out}$, the other is for large $k_{out}$ in the tail.}.

In our collected population of inter-connected blogs, the maximum
in-degree is $13,342$, whereas a majority of blogs just have a few
incoming links (see Tab.~\ref{tab:degree}). The power-law
distribution of in-degree indicates that many common bloggers
preferentially add links to their favorite celebrities' blogs and
such preferential behavior results in the power-law distribution
of the in-degree just as the BA model describes. We found that a
significant fraction of the blogs, that is 32.6\%, have no
outgoing links to other blogs. Further, considerable fraction of
blogs have only a few outgoing connections (see
Tab.~\ref{tab:degree}). That's to say, most of the bloggers are
unknown to public and they are not active enough in the
blogosphere (have small or null outgoing links).

The average degree $\langle k \rangle$ of such blogging network is
$9.0$, that's to say, for each node in such social networks has an
average of $9$ neighbors. Furthermore, the average in and out
degrees $\langle k_{in} \rangle = \langle k_{out} \rangle = 4.5$.
Although there are millions of connections presented in the social
network, as aforementioned, about $28.7\%$ of them are symmetric
and most of the symmetric links are between the blogs of bloggers
who get acquainted with each other in the blogosphere. So this
proves that such blogging network is asymmetric one: while a node
tends to link to a famous node, it is seldom the case that the
famous node would link to this node either.

\begin{table}
\caption{\label{tab:degree}Percentage of blogs with null, 1, 2 and
3 out and in degrees. Note that a large fraction of blogs have
only small in and out degrees. Since our blogging network was
crawled along the directed links, the in-degree of blogs is at
least 1. }
\begin{ruledtabular}
\begin{tabular}{ccccc}
k= & 0 & 1 & 2 & 3\\
\hline
In & 0 & 48.4\% & 18.1\% & 9.8\%\\
Out & 32.6\% & 14.2\% & 9.7\% & 7.5\%\\

\end{tabular}
\end{ruledtabular}
\end{table}

\begin{table}
\caption{\label{tab:dc}Correlation coefficients for the degrees at
either side of an edge. Negative figures indicate that poorly
connected nodes tend to link to highly connected nodes while
positive values suggest that nodes with even connectivity are
likely to connect to each other.}
\begin{ruledtabular}
\begin{tabular}{ccccc}
$r$ & $r_{in-in}$ & $r_{in-out}$ & $r_{out-in}$ & $r_{out-out}$\\
\hline
$-$0.497 & $-$0.035 & 0.041 & $-$0.034 & 0.113\\
\end{tabular}
\end{ruledtabular}
\end{table}

Table~\ref{tab:dc} displays the correlation coefficients of
different types of degree-degree correlations for the crawled down
blogging network. Correlations are measured by the Pearson's
correlation coefficient $r$ for the degrees at either side of an
edge as suggested by Mark Newman \cite{degree-cor}.
\begin{equation}
r=\frac{\langle k_{to}k_{from} \rangle - \langle k_{to}
\rangle\langle k_{from} \rangle}{\sqrt{\langle k_{to}^2 \rangle -
\langle k_{to} \rangle^2}\sqrt{\langle k_{from}^2 \rangle -
\langle k_{from} \rangle^2}}
\end{equation}
where $k_{to},k_{from}$ could be four possible combinations of in
and out degrees of an edge.

Networks with assortative mixing pattern are those in which nodes
with large degree tend to be connected to other nodes with many
connections and vice visa. Technical and biological networks are
in general disassortative, while social networks are often
assortatively mixed as demonstrated by the study on scientific
collaboration networks \cite{degree-cor}. Blogging network,
however, presents disasortative mixing pattern when directions are
not considered. Positive mixing are shown for $r_{in-out}$ and
$r_{out-out}$ in our case. Positive $r_{in-out}$ means active
bloggers in the community (have large $k_{out}$) tend to associate
with those who succeed in promoting themselves in the community
(have high $k_{in}$), while a large $r_{out-out}$ suggests that
the active bloggers preferentially link to each other. Internet
dating community, a kind of social networks embedded in a
technical one, and peer to peer (P2P) social networks are similar
to our case, displaying a significant disassortative mixing
pattern \cite{Internet_dating,P2P}.

\begin{figure}
\includegraphics[width=8cm]{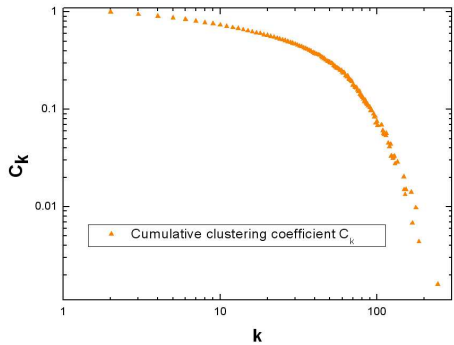}
\caption{\label{fig:clustering}(Color online) Cumulative
distribution of clustering coefficient of blogs.}
\end{figure}

The length of average shortest path $\langle l \rangle$ is
calculated, which is the mean of geodesic distance between any
pairs that have at least a path connecting them. In this case,
$\langle l \rangle = 6.84$. That means on average one only needs
to click $7$ times from one blog site to any other blog site in
the blogosphere. And the diameter $D$ of this social networks
which is defined as the maximum of the shortest path length, is
$27$. Because such blogging network is directed, the clustering
coefficient is not easy to be computed. One way to avoid this
difficulty is to make the network undirected. Firstly, the one-way
connections were removed from the network; secondly, the isolated
nodes were deleted from the graph. By doing so, the bidirectional
graph with $122,470$ nodes was obtained to compute the clustering
coefficient. The mean degree of this undirected networks
$k_{undirected}$ is 3.28. According to the definition of
clustering coefficient in undirected network,
$C_i=\frac{2E_i}{k_i(k_i-1)}$, that is the ratio between the
number $E_i$ of edges that actually exits between these $k_i$
neighbor nodes of node $i$ and the total number $k_i(k_i-1)$. The
clustering coefficient of the whole network is the average of all
individual $C_i$'s. We found the clustering coefficient
$C=0.1490$, order of magnitude much higher than that of a
corresponding random graph of the same size
$C_{rand}=3.28/122470=0.0000268$. Besides, the degree-dependent
local clustering coefficient $C(k)$ is averaging $C_i$ over
vertices of degree $k$. Fig.~\ref{fig:clustering} plots the
cumulative distribution of $C(k)$ from the undirected blogging
network. However, it is hard to declare a clear power law in our
case. Nevertheless, the nonflat clustering coefficient
distributions shown in the figure suggests that the dependency of
$C$ on $k$ is nontrivial, and thus points to some degree of
hierarchy in the networks. Consequently, it is demonstrated that
the average shortest path length is far smaller than the logarithm
of the network size in such blogging network. In addition, the
network has relatively high clustering coefficient. Thence, the
blogging network of inter-connected blogs has small-world effect.
This small-world phenomenon is also consistent with the former
small-world discovery about the WWW.

\begin{figure}
\includegraphics[width=8cm]{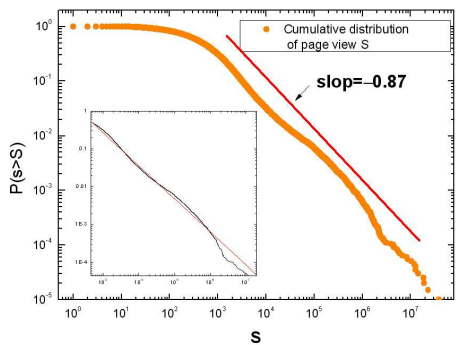}
\caption{\label{fig:pv}(Color online) Cumulative distribution of
the fraction of blogs of which the number of page views is more
than $S$. The stright-line is of slope $-0.87$ for comparison with
the distribution.}
\end{figure}

To evaluate the popularity of the blogs, the cumulative
distribution of the number of page views of blog sites is figured
out (see Fig.~\ref{fig:pv}). For small page view $S$ ($S\le 500$),
there exists saturation. However, for large $S$, the fraction of
blogs that have more than total $S$ page views obeys a power-law
form as $P(s>S)\sim S^{-\tau}$, where $\tau=0.87\pm
6.56\times10^{-4}$. Immediately, one could get the distribution of
page views of blogs as $p(s)\sim s^{\mu}$, where
$\mu=\tau+1=1.87$. In our case of $200339$ nodes, only ten blogs's
page views exceed tens of millions, while most of the remanent
have only tens of thousands page views. This heavy-tailed
distribution indicates that most of the readers are attracted by
the celebrated bloggers and contribute page views to their blogs.
However, minority of the grassroots' blogs could gain public
attention in the blogosphere. In this sense, some kind of
inequality develops: the richer gets richer while the poorer gets
poorer. Thus, social technologies not only enhance the
communications between distant people, but also facilitate the
inequality between the celebrated and the commons. From this
respect, the blogosphere might be a good paradigm for studying the
emergence of such inequality.

\subsection{Conclusion remarks and future work}
In summary, the sub-ecosystem of global blogosphere is scrutinized
to reveal the underlying link pattern and the popularity of the
blogs. We found that the blogging community has small-world
property. In addition, the in-degree and out-degree distributions
follow power-law forms. Calculations on degree-degree correlations
show that blogging networks are in general disassortative mixing,
except that active bloggers are connected between each other and
by the ones with high in-degree. The fraction of number of page
views of blogs also obeys a power law. Although our crawled down
blogging network is static whereas the nature of blogosphere is
the dynamical and evolving one, our observations and statistical
analysis might be the first step to such ecosystem. However, what
has been done is not enough. There are still various aspects of
blogoshpere to be investigated. Recently, a new technique called
collaborative tagging gains ground in blogging community because
it could steer bloggers to effectively share tremendous amounts of
information and find the useful information
\footnote{http://en.wikipedia.org/wiki/Tags}. It is of some merit
to study tag co-occurrence to reveal the universal characteristics
of users' tagging behavior \cite{tag}. Moreover, the fascinating
phenomenon of arising hot discussion topics is worth examining to
dig out the intrinsic features of collective behaviors in
blogosphere. And also the recommendations rules of blogging
creates powerful collaborative filtering. Thus blogosphere would
be the suitable one to study collaborative filtering effect.
Additionally, detecting the latent community structures in
blogosphere would be meaningful. And also, it is interesting to
study information, like rumors, propagation in this ecosystem. In
short, self-organized blogosphere is a good paradigm for
understanding varieties of facets of behavior pattern of bloggers
in such ecosystem.

\begin{acknowledgments}
The authors are partly supported by National Natural Science
Foundation of China (NNSFC) under Grant No.10372002 and
No.60528007 and National 973 Program under Grant No.2002CB312200.
\end{acknowledgments}

\end{document}